\documentstyle[12pt]{article}
\title{Quasi-triangular structures on \\
Hopf algebras with positive bases}
\begin{document}
\author{Jiang-Hua Lu\thanks{Research is partially supported by an
NSF postdoctorial fellowship and NSF grant DMS 9803624.}
\\
Department of Mathematics\\ 
University of Arizona\\
\\
Min Yan, \hspace{5mm}
Yong-Chang Zhu\thanks{Research is supported by UGC Earmark Grant
HKUST 629/95P. }
  \\
Department of Mathematics\\ 
Hong Kong University of Science and Technology} 
\date{}
\maketitle

\newcommand{\la}{\leftarrow}
\newcommand{\ra}{\rightarrow}
\newcommand{\da}{\downarrow}
\newcommand{\ua}{\uparrow}
\newcommand{\uda}{\updownarrow}
\newcommand{\sub}{\subset}
\newcommand{\lra}{\longrightarrow}
\newcommand{\lla}{\longleftarrow}
\newcommand{\La}{\Leftarrow}
\newcommand{\Ra}{\Rightarrow}
\newcommand{\Ua}{\Uparrow}
\newcommand{\Da}{\Downarrow}
\newcommand{\Uda}{\Dpdownarrow}
\newcommand{\Lla}{\Longleftarrow}
\newcommand{\Lra}{\Longrightarrow}
\newcommand{\Llra}{\Longleftrightarrow}
\newcommand{\pa}{\partial}
\newcommand{\lcpr}{ \!>\!\!\!\triangleleft\, }
\newcommand{\rcpr}{ \,\triangleright\!\!\!<\! }

\newtheorem{thm}{Theorem}
\newtheorem{lem}{Lemma}
\newtheorem{cor}[lem]{Corollary}
\newtheorem{prop}[lem]{Proposition}
\newtheorem{df}{Definition}
\newtheorem{rmk}{Remark}
\newtheorem{eg}{Example}

\tableofcontents

\newpage

\begin{abstract}
A basis $B$ of a finite dimensional Hopf algebra $H$ is said to
be  positive if all the structure constants  of $H$ relative to $B$
are  non-negative. A quasi triangular structure 
$R\in H\otimes H$ is said to be positive with respect to $B$ if it has
non-negative coefficients in the basis $B \otimes B$ of 
$H\otimes H$. In our earlier work, we have classified all finite
dimensional Hopf
algebras with positive bases. In this paper, we 
classify positive quasi-triangular structures on such Hopf
algebras. A consequence of this
classification is a new way of constructing set-theoretical
solutions of the Yang-Baxter equation.
\end{abstract}

\section{Introduction}

Consider a finite dimensional Hopf algebra $H$ with a basis
such that all the structure constants with respect to this
basis  are non-negative. We call such a basis 
{\em  positive}. In \cite{lyz:positive} we showed that any such
Hopf algebra is isomorphic to the bicrossproduct Hopf algebra
$H(G;G_+,G_-)$ coming from a 
unique factorization $G=G_+G_-$ of a finite group $G$ (see
Section \ref{sec_gd} for the definition of $H(G;G_+,G_-)$).
We also showed 
that such Hopf
algebras are exactly the linearizations of Hopf algebras in the
category of sets with correspondences as morphisms.

In this paper, we classify all quasi-triangular
structures $R \in H \otimes H$ on the Hopf algebra 
$H=H(G;G_+,G_-)$ that are {\it positive}
 in the sense that the coefficients of $R$
in the basis $G \otimes G$ of $H\otimes H$ are non-negative. 
We also show that they
are all quasi-equivalent to certain normal forms. 
It turns out that such an $R$ is 
set-theoretical, in the sense that it is the linearization
of a bisection ${\cal R}$ of a product groupoid 
$\Gamma\times\Gamma$ satisfying the {\em groupoid-theoretical}
Yang-Baxter equation. Consequently, for every
$\Gamma$-set $X$, ${\cal R}$ induces a set-theoretical
solution of the Yang-Baxter equation on $X$. This fact
motivates our work in
\cite{lyz:setybe}, where we give a general way of constructing
set-theoretical solutions of the  Yang-Baxter equation that 
includes as special cases the earlier ones by Weinstein and
Xu  in \cite{w-x:R} and by Etingof, Schedler and Soloviev  in
\cite{e-t-s:set}. We also classify positive triangular structures
on $H(G; G_+, G_-)$ and recover a construction of such
structures by Etingof and Gelaki \cite{eg}.

{\bf Acknowledgment} We would like to thank Professor P. Etingof
for helpful comments.

\section{Positive quasi-triangular structures} 
\label{sec_gd}

A {\em unique factorization} $G=G_+G_-$ of a group $G$
consists of two subgroups $G_+$ and $G_-$ such that any $g\in
G$ can be written as $g=g_+g_-$ for unique $g_+\in G_+$ and 
$g_-\in G_-$. By considering the inverse map, we see that
for every $g\in G$, there are unique
$g_+,\bar{g}_+\in G_+$ and $g_-,\bar{g}_-\in G_-$ such that
\[
g=g_+g_-=\bar{g}_-\bar{g}_+.
\]
We will denote $(g_+)^{-1} \in G_+$ and $(g_-)^{-1} \in G_-$ simply by
$g_+^{-1}$ and $g_-^{-1}$.

The unique factorization induces the following
actions of $G_+$ and $G_-$ on each other (from left and from right) 
\begin{eqnarray*}
G_-\times G_+\rightarrow  G_+, & &
(\bar{g}_-,\bar{g}_+) \mapsto 
g_+=\;^{\bar{g}_-}\bar{g}_+, \\
G_-\times G_+\rightarrow G_-, & &
(\bar{g}_-,\bar{g}_+) \mapsto 
g_-=\;\bar{g}_-^{\,\,\bar{g}_+}, \\
G_+\times G_-\rightarrow G_+, & &
(g_+,g_-) \mapsto \bar{g}_+=\;g_+^{\,\,g_-}, \\
G_+\times G_-\rightarrow G_-, & &
(g_+,g_-) \mapsto \bar{g}_-=\;^{g_+}g_-.
\end{eqnarray*}
By definition, we have
$g_+g_-=\left(^{g_+}g_-\right)\left(g_+^{\,\,g_-}\right)$
and 
$g_-g_+=\left(^{g_-}g_+\right)\left(g_-^{\,\,g_+}\right)$.
Moreover, the actions have the following properties
\begin{equation}\label{compatible}
\left\{
\begin{array}{ll}
\;^{g_+}(g_-h_-)=\;^{g_+}g_-\;\;^{(g_+^{\;g_-})}h_-, &
(h_+g_+)^{g_-}=h_+^{(\;^{g_+}g_-)}\;g_+^{\;g_-}, \\
\;^{g_-}(g_+h_+)=\;^{g_-}g_+\;\;^{(g_-^{\;g_+})}h_+, &
(h_-g_-)^{g_+}=h_-^{(\;^{g_-}g_+)}\;g_-^{\;g_+}.
\end{array}
\right.
\end{equation}
\begin{equation}\label{inverse}
\left\{
\begin{array}{ll}
\left(g_+\,^{g_-}\right)^{-1}=\,^{(g_-^{-1})}g_+^{-1}, &
\left(\,^{g_-}g_+\right)^{-1}=g_+^{-1}\,^{(g_-^{-1})}, \\
\left(g_-\,^{g_+}\right)^{-1}=\,^{(g_+^{-1})}g_-^{-1}, &
\left(\,^{g_+}g_-\right)^{-1}=g_-^{-1}\,^{(g_+^{-1})}.
\end{array}
\right.
\end{equation}

A Hopf algebra $H(G;G_+,G_-)$ can be constructed from a unique
factorization $G=G_+G_-$ of a finite group  \cite{mj:book}
\cite{ta:bismash}. More precisely, $H(G;G_+,G_-)$ is the vector space 
spanned by the set $G$ together with the 
following Hopf algebra structure (we use $\{g\}$ to
denote the group element $g\in G$ considered as an element of 
$H(G;G_+,G_-)$):
\[
\left\{\begin{array}{ll}
\mbox{multiplication:} & 
\{g\} \{h\} = \delta_{g_+^{g_-}, h_+}\{gh_-\} \\
\noalign{\vskip 3pt}
\mbox{unit:} & 
1 = \sum_{g_+ \in G_+} \{g_+\} \\
\noalign{\vskip 5pt}
\mbox{co-multiplication:} &
\Delta\{g\} = \sum_{h_+ \in G_+}
\{g_+ h_+^{-1} (^{h_+}g_-)\} \otimes \{h_+g_-\} \\
\noalign{\vskip 5pt}
\mbox{co-unit:} &
\epsilon\{g\} = \delta_{g_+, e} \\
\noalign{\vskip 5pt}
\mbox{antipode:} & S\{g\} = \{g^{-1}\}
\end{array}\right.
\]
We remark that the algebra structure on $H(G; G_+, G_-)$ is that of
the cross-product of the group algebra ${\bf C} G_{-}$ of
$G_{-}$ and the function algebra ${\bf C}(G_+)$ of $G_+$
with respect to the above right action of $G_-$ on $G_+$, and similarly
for its co-algebra structure.
The Hopf algebra $H(G;G_+,G_-)$ has $G$ as the obvious
positive basis. In \cite{lyz:positive}, we proved that all
finite dimensional Hopf algebras with positive bases are
of the form $H(G;G_+,G_-)$.

Recall that a {\em quasi-triangular structure} on a Hopf
algebra $H$ is an invertible element $R \in H \otimes H$ such
that
\begin{eqnarray}
\label{eq_prime}
& & \tau\Delta(g) \,  = \,  R \Delta(g) R^{-1}, \qquad  
\mbox{for all $g \in H$} \\
\label{eq_delta-R}
& & (\Delta \otimes id)R \, = \, R_{13}R_{23}, \qquad 
(id \otimes \Delta) R \,  = \,  R_{13}R_{12}\\
\label{eq_epsilon-R}
& & (\varepsilon \otimes id) R \,  = \,  
(id \otimes \varepsilon)R \, = \,  1,
\end{eqnarray}
where $\tau(a\otimes b)=b\otimes a$. 
In the special case $H=H(G;G_+,G_-)$, we say that an 
element $R\in H\otimes H$ is {\em positive} if
$R$ is a non-negative linear combination of the basis elements 
$\{g\} \otimes \{h \}$, $g, h \in G$.

\begin{thm}\label{r-classification}
Let $G= G_+ G_-$ be a unique factorization of a finite group
$G$. Let $\xi,\eta:G_+\rightarrow G_-$ be two group
homomorphisms,  and denote  
\[
G_+'=\{u\xi(u^{-1}): \, u \in  G_+\},
\qquad
G_+''=\{\eta(u^{-1})u: \, u\in G_+\}, 
\]
\[
F(u\xi(u^{-1}))=\eta(u)u^{-1}: \quad G_+'\rightarrow G_+''.
\]
Suppose that the following conditions are satisfied,

{\em (a)} Both $G_+'$ and $G_+''$ are normal subgroups of $G$;

{\em (b)} $F$ is a group isomorphism.

\noindent Then 
\[
R=\sum_{u,v\in G_+}
  \{u\left(\eta(v)^{u}\right)^{-1}\}\otimes \{v\xi(u)\}
\]
is a positive quasi-triangular structure. Conversely, 
every positive quasi-triangular structure on $H(G;G_+,G_-)$ is
given by the construction above.
\end{thm}

The following proposition rephrases the conditions (a) and
(b) in a more concrete form. Theorem \ref{r-classification} is
a consequence of this proposition and Propositions 
\ref{positiveR} - \ref{lem4} of Section \ref{proof}.

\begin{prop}\label{condition-construction1}
Let $G=G_+G_-$ be a unique factorization. 
Let $\xi,\eta: G_+\rightarrow G_-$ be two group homomorphisms.
Then the conditions {\em (a)} and {\em (b)} in Theorem
\ref{r-classification} are equivalent to 
\begin{eqnarray}
\xi(u)^v & = & \xi(u^{\eta(v)}), \label{etaright} \\
\,^u\eta(v) & = & \eta(^{\xi(u)}v), \label{xileft} \\
uv & = & (^{\xi(u)}v)(u^{\eta(v)}), \label{etaxi}  \\
\xi(^x u)x\,^u & = & x\xi(u), \label{etaneg'} \\
\eta(^x u)x\,^u & = & x\eta(u), \label{xineg'}  
\end{eqnarray}
for all $u,v\in G_+$, and $x\in G_-$. Moreover, each of these properties
is also equivalent to the corresponding property below
\begin{eqnarray}
\,^v\xi(u) & = & \xi(^{\eta(v)}u), \label{etaleft} \\
\eta(v)^u & = & \eta(v^{\xi(u)}), \label{xiright} \\
uv & = & (^{\eta(u)}v)(u^{\xi(v)}), \label{xieta}  \\
\,^u x\xi(u\,^x) & = & \xi(u)x, \label{etaneg} \\
\,^u x\eta(u\,^x) & = & \eta(u)x. \label{xineg} 
\end{eqnarray}
\end{prop}

\noindent{\em Proof}: We first prove that (a) and (b) imply 
$(\ref{etaright}-\ref{xineg'})$.

Since $G_+'$ is normal, for any
$u\in G_+$ and $x\in G_-$, we can find $v\in G_+$ such that
$xu\xi(u^{-1})=v\xi(v^{-1})x$. By the unique factorization
$G=G_+G_-$, we have
\[
\,^xu=v,\qquad (x^u)\xi(u^{-1})=\xi(v^{-1})x.
\]
This implies (\ref{etaneg'}). Similarly, the fact that $G_+''$
is normal implies (\ref{xineg'}).

Since 
\begin{eqnarray}
u\xi(u^{-1})v\xi(v^{-1})
& = &  
u(^{\xi(u^{-1})}v)(\xi(u^{-1})^{v})\xi(v^{-1}), \label{rel52} \\
\eta(u)u^{-1}\eta(v)v^{-1}
& = & 
\eta(u)(^{u^{-1}}\eta(v))((u^{-1})^{\eta(v)})v^{-1},
\label{rel51}  
\end{eqnarray}
and $F$ is a homomorphism, we have 
$(u(^{\xi(u^{-1})}v))^{-1}=((u^{-1})^{\eta(v)})v^{-1}$. This is
exactly (\ref{etaxi}).

By (\ref{rel52}) and the fact that $G_+'$ is a subgroup, we have
\begin{equation}\label{rel54}
\xi(u(^{\xi(u^{-1})}v))^{-1}
=(\xi(u^{-1})^{v})\xi(v^{-1}).
\end{equation}
Since $u(^{\xi(u^{-1})}v)\stackrel{(\ref{etaxi})}{=}
v((u^{-1})^{\eta(v)})^{-1}$,
the left side of (\ref{rel54}) is 
$\xi((u^{-1})^{\eta(v)})\xi(v^{-1})$.
Therefore (\ref{rel54}) is equivalent to 
$\xi((u^{-1})^{\eta(v)})=\xi(u^{-1})^{v}$, which is
exactly (\ref{etaright}).
Similarly, the fact that $G_+''$ is a subgroup implies
(\ref{xileft}).

The converse that $(\ref{etaright}-\ref{xineg'})$ implies (a) and
(b) can be proved similarly, by making use of the same
computations. The conditions (\ref{etaright}), (\ref{xileft}),
and (\ref{etaxi}) imply that $G_+'$ and
$G_+''$ are subgroups. The conditions (\ref{etaneg'}) and
(\ref{xineg'}) imply that $G_+'$ and $G_+''$ are normal subgroups.
The condition (\ref{etaxi}) implies that $F$ is a group
isomorphism.

Finally, if (\ref{etaright}) holds, then we have
\[
\,^v\xi(u)
\stackrel{(\ref{inverse})}{=}
\left(\xi(u^{-1})^{v^{-1}}\right)^{-1}
\stackrel{(\ref{etaright})}{=}
\xi\left((u^{-1})^{\eta(v^{-1})}\right)^{-1}
\stackrel{(\ref{inverse})}{=}
\xi(^{\eta(v)}u),
\]
which is (\ref{etaleft}). The same idea shows that each of
$(\ref{xiright}-\ref{xineg})$ is equivalent to the corresponding
property in $(\ref{xileft}-\ref{xineg'})$.

\hfill$\Box$

Next we give an alternative description for the data
$(G=G_+G_-,\xi,\eta)$ used in Theorem 
\ref{r-classification}.

Let $G_-$ be a group acting on another group $A$ as
automorphisms, with the action denoted by $(x,a)\mapsto x\cdot a:
G_-\times A\rightarrow A$. Then we have the semi-direct product
group $G= A\lcpr G_-$, with the group structure given by
\[
(ax)(by)=a(x\cdot b)xy,\qquad a,b\in A,\quad x,y\in G_-.
\]
A map $\zeta: A\rightarrow G_-$ is called  a {\em 1-cycle} if
\begin{equation}\label{cycle-condition}
\zeta(a)\zeta(b)=\zeta(a(\zeta(a)\cdot b)).
\end{equation}
The 1-cycle condition (\ref{cycle-condition}) is equivalent to the
fact that $\{a\zeta(a):\,a\in A\}$ is a subgroup of $G$. Moreover, if
$\zeta$ is bijective, then $\zeta^{-1}: G_-\rightarrow A$ is a
$1$-cocycle of $G_-$ with coefficient in $A$ as defined in
\cite{eg}.

\begin{thm}\label{r-cycle}
There is a one-to-one correspondence between triples
$(G=G_+G_-,\xi,\eta)$ satisfying the conditions of Theorem 
\ref{r-classification} and the triples 
$(G=A\lcpr G_-,\zeta,F)$, where $\zeta$ is a cycle of $G_-$ with
coefficients in $A$ and $F$ is an automorphism of $G$,
satisfying

{\em (a)} $F(x)=x$ for any $x\in G_-$;

{\em (b)} $F(a)a\in G_-$ for any $a\in A$.
\end{thm}

Specifically, the correspondence is the following.
Given $(G=G_+G_-,\xi,\eta)$, we define
\begin{equation}\label{1to2}
\begin{array}{l}
A=G_+'=\{u\xi(u^{-1}):\, u\in G_+\}; \\
\noalign{\vskip 5pt}
F(u\xi(u^{-1})x)=\eta(u)u^{-1}x, \quad
\mbox{for $u\in G_+$ and $x\in G_-$}; \\
\noalign{\vskip 5pt}
\zeta(u\xi(u^{-1}))=\xi(u), \quad
\mbox{for $u\in G_+$}.
\end{array}
\end{equation}
Moreover, since $A$ is a normal subgroup, conjugations by elements in
$G_-$ give an action of $G_-$ on $A$ as automorphisms.
Conversely, given $(G=A\lcpr G_-,\zeta,F)$, we define
\begin{equation}\label{2to1}
\begin{array}{l}
G_+=\{a\zeta(a):\, a\in A\}; \\
\noalign{\vskip 5pt}
\xi=P|_{G_+}; \\
\noalign{\vskip 5pt}
\eta=P\circ F^{-1}|_{G_+},
\end{array}
\end{equation}
where $P$ is the natural homomorphism 
$G=A\lcpr G_-\rightarrow G_-$,
$ax\mapsto x$.

\bigskip

\noindent{\em Proof of Theorem \ref{r-cycle}}: First we show that
if $(G=G_+G_-,\xi,\eta)$ satisfies the conditions of Theorem 
\ref{r-classification}, then the construction (\ref{1to2}) is as
described in Theorem \ref{r-cycle}. 

For $a=u\xi(u^{-1})$, we have $a\zeta(a)=u$. Therefore the subset
$\{a\zeta(a):\,a\in A\}=G_+$ is a subgroup of $G$. This implies $\zeta$
is a 1-cycle.

By its very definition, $F$ is a homomorphism if and only if
$F:G_+'\rightarrow G_+''$ is an equivariant map with respect to
the $G_-$-actions defined by conjugations. For $x\in G_-$
and $a=u\xi(u^{-1})\in G_+'$, the action of $x$ on $a$
is 
\[
x\cdot a = xu\xi(u^{-1})x^{-1} = v\xi(v^{-1}),
\qquad v=(xu)_+=\,^xu,
\]
and the action of $x$ on $F(a)$ is 
\[
x\cdot F(a)=x\eta(u)u^{-1}x^{-1}=\eta(w)w^{-1},
\qquad w=(xu)_+=\,^xu.
\]
We conclude from this that $v=w$ and $F(x\cdot a)=x\cdot F(a)$.

Finally, the condition (a) follows from the definition, and
(b) follows from
\[
F(u\xi(u^{-1})) u\xi(u^{-1})=\eta(u)u^{-1}u\xi(u^{-1})
=\eta(u)\xi(u^{-1}).
\]

Now we turn to the construction (\ref{2to1}).

First of all, since $\zeta$ is a 1-cycle, we know
$G_+$ is a subgroup of $G$. Moreover, for any $a\in A$ and $x\in G_-$,
the decomposition $ax=(a\zeta(a))(\zeta(a)^{-1}x)$
gives the unique factorization $G=G_+G_-$.

Since $P$ and $F^{-1}$ are homomorphisms, $\xi$ and $\eta$ are
also homomorphisms. 

We express an element in $G$ as $ux$ for unique $u\in G_+$ and
$x\in G_-$. The element is in $A$ if and only if it is in the
kernel of $P$. Since $P(ux)=P(u)P(x)=\xi(u)x$, we see that
$A$ consists of elements of the form $u\xi(u^{-1})$, $u\in G_+$.
In other words, we have $A=G_+'$, which in particular implies
$G_+'$ is a normal subgroup. Similarly, by considering those elements
$xu$ in the kernel of
$P\circ F^{-1}$, we conclude that $G_+''=F(A)$. Since $F$ is an
automorphism, $G_+''$ is also a normal subgroup.

Since $F(G_+')=G_+''$, for any $u\in G_+$, we can find $v\in G_+$ such
that $F(u\xi(u^{-1}))=\eta(v)v^{-1}$. Then by condition (b), we
have 
$\eta(v^{-1})vu\xi(u^{-1})=F(u\xi(u^{-1}))u\xi(u^{-1})\in G_-$.
This implies $uv\in G_-$. On the other hand, $u,v\in G_+$
implies $uv\in G_+$. Therefore by the unique factorization
$G=G_+G_-$, we have $uv=e$. Consequently, the formula 
$F(u\xi(u^{-1}))=\eta(u)u^{-1}$ holds.

\hfill$\Box$

\section{Proof of the classification theorem}
\label{proof}

We prove Theorem \ref{r-classification} in this section.

\begin{prop}\label{positiveR}
Suppose that 
$R\in H(G;G_+,G_-)\otimes H(G;G_+,G_-)$ is invertible and positive, and
suppose that $R^{-1}$ is also positive. Then there is a subset 
${\cal R}\subset G\times G$ and a positive valued function 
$r:{\cal R}\rightarrow {\bf R}^{>0}$ such that
\begin{enumerate}
\item The restriction of the map $(g,h)\rightarrow (g_+,h_+): 
      G\times G\rightarrow G_+\times G_+$ to ${\cal R}$ is a bijection;
\item $R=\sum_{(g,h)\in {\cal R}}r(g,h)\{g\}\otimes \{h\}$.
\end{enumerate}
\end{prop}

\noindent{\em Proof}: 
The positivity assumption implies that
if the coefficients of $\{g\}\otimes\{h\}$ in $R$ and 
of $\{k\}\otimes \{l\}$ in $R^{-1}$ are non-zero, then either the
multiplicability conditions $g_+^{g_-}=k_+$ and $h_+^{h_-}=l_+$ are not
satisfied, or the coefficient of the product
$(\{g\}\otimes\{h\})(\{k\}\otimes\{l\})=\{gk_-\}\otimes\{hl_-\}$ in 
$RR^{-1}=e\otimes e=\sum_{u,v\in G_+}u\otimes v$ is non-zero. Similarly
for $R^{-1}R$. The proposition is then a consequence of these two
facts.

Instead of going through the details of the argument, we note that the
proposition is a consequence of a general fact about bisections of
groupoids. See the discussion after the proof of Proposition \ref{bisect}
for such a conceptual proof of the proposition.

\hfill$\Box$

\begin{prop}\label{lem1}
Suppose that $\xi,\eta:G_+\rightarrow G_-$ are two group homomorphisms.
Suppose also that $r:G_+\times G_+\rightarrow {\bf R}^{>0}$ is a function
such that for any $u,v,w\in G_+$, the equalities
$(\ref{etaleft})$, $(\ref{xiright})$, and the following are
satisfied
\begin{eqnarray}
r(uw,v) & = & r(u,v)r(w,v^{\xi(u)}), \label{reta} \\
r(u,wv) & = & r(u,v)r(^{\eta(v)}u,w). \label{rxi} 
\end{eqnarray}
Then
\begin{equation}\label{r-form'}
R=\sum_{u,v\in G_+}r(u,v)
\{u\left(\eta(v)^{u}\right)^{-1}\}\otimes \{v\xi(u)\}
\end{equation}
satisfies $(\ref{eq_delta-R})$. Conversely, if $R$ is invertible,
positive, satisfies $(\ref{eq_delta-R})$, and $R^{-1}$ is
also positive, then $R$ is given by the
construction above.
\end{prop}

\noindent{\em Proof}: Suppose that $R$ is invertible,
positive, and $R^{-1}$ is also positive. Then Proposition
\ref{positiveR} implies that 
\begin{equation}\label{r-form}
R=\sum_{u,v\in G_+}r(u,v)\{u\phi(u,v)\}\otimes \{v\psi(u,v)\},
\end{equation}
where $\phi, \psi: G_+\rightarrow G_-$ are two maps, and
$r:G_+\times G_+\rightarrow {\bf R}^{>0}$ is a positively valued
function. From (\ref{r-form}), we have
\[
(\Delta\otimes id)R = \sum_{u,v,w\in G_+}
r(u,v)
\{uw^{-1}(^w \phi(u,v))\}\otimes
\{w\phi(u,v)\}\otimes
\{v\psi(u,v)\}
\]
and
\begin{eqnarray*}
R_{13}R_{23}  & \displaystyle = \sum_{u,v,w\in G_+} &
r(u,v)r(w,v^{\psi(u,v)})
\{u\phi(u,v)\}\otimes 
\{w\phi(w,v^{\psi(u,v)})\}  \\
&&
\otimes
\{v\psi(u,v)\psi(w,v^{\psi(u,v)})\}.
\end{eqnarray*}
Then it is easy to see that $(\Delta\otimes id)R=R_{13}R_{23}$ means
\begin{eqnarray}
r(uw,v) & = & r(u,v)r(w,v^{\psi(u,v)}) \label{rel1} \\
\,^{w} \phi(uw,v) & = & \phi(u,v) \label{rel2} \\
\phi(uw,v) & = & \phi(w,v^{\psi(u,v)}) \label{rel3} \\
\psi(uw,v) & = & \psi(u,v)\psi(w,v^{\psi(u,v)}) \label{rel4}
\end{eqnarray}
for all $u,v,w\in G_+$.
Similarly, we see that $(id\otimes\Delta)R=R_{13}R_{12}$ means
\begin{eqnarray}
r(u,wv) & = & r(u,v)r(u^{\phi(u,v)},w) \label{rel5} \\
\phi(u,wv) & = & \phi(u,v)\phi(u^{\phi(u,v)},w) \label{rel6} \\
\,^{v}\psi(u,wv) & = & \psi(u^{\phi(u,v)},w) \label{rel7} \\
\psi(u,wv) & = & \psi(u,v) \label{rel8}
\end{eqnarray}
for all $u,v,w\in G_+$.

Equation (\ref{rel8}) implies that 
\begin{equation}\label{rel9}
\psi(u,v)=\xi(u)
\end{equation}
for a map $\xi: G_+\rightarrow G_-$. Then (\ref{rel4}) becomes 
$\xi(uw)=\xi(u)\xi(w)$, i.e., $\xi$ is a group
homomorphism.

Equation (\ref{rel2}) implies that $\phi(u,v)=\,^{u^{-1}}\phi(e,v)$.
Therefore we introduce $\eta(v)=\phi(e,v)^{-1}: G_+\rightarrow G_-$ and
have
\begin{equation}\label{rel10}
\phi(u,v)
=\,^{u^{-1}}(\eta(v)^{-1})
\stackrel{(\ref{inverse})}{=}\left(\eta(v)^{u}\right)^{-1}.
\end{equation}
Moreover, we have
\begin{equation}\label{rel11}
u^{\phi(u,v)}
=
u^{(^{u^{-1}}(\eta(v)^{-1}))}
\stackrel{(\ref{compatible})}{=}
\,^{\eta(v)}u.
\end{equation}
Then by (\ref{rel10}) and (\ref{rel11}), equation
(\ref{rel6}) becomes
\[
\eta(wv)^u 
=\eta(w)^{(^{\eta(v)}u)}
 \eta(v)^u.
\]
Taking $u=e$, we see that $\eta$ is a group homomorphism. By
making use of this fact, the equation above becomes (\ref{compatible}),
which is always satisfied.

By (\ref{rel9}) and (\ref{rel11}), equation
(\ref{rel7}) becomes (\ref{etaleft}).
By (\ref{rel10}), equation (\ref{rel3}) becomes
$\eta(v)^{uw}=\eta(v^{\xi(u)})^w$. Applying the right action by $w^{-1}$,
we have (\ref{xiright}).
Finally, by (\ref{rel11}), equations (\ref{rel1}) and
(\ref{rel5}) become (\ref{reta}) and (\ref{rxi}).

\hfill$\Box$

\begin{prop}\label{lem2}
Suppose that $\xi$ and $\eta:G_+\rightarrow G_-$ 
are two group homomorphisms. Suppose also that
$r:G_+\times G_+\rightarrow {\bf R}^{>0}$ is a function such
that for any $u,v,w\in G_+$ and $x\in G_-$, the equalities
$(\ref{etaright}-\ref{xineg'})$ and the following are satisfied
\begin{equation}\label{rneg}
r(u,v) = r(u^x,v^{(^u x)}). 
\end{equation}
Then $(\ref{r-form'})$ is a positive quasi-triangular structure
on $H(G;G_+,G_-)$. 
Conversely,  any positive quasi-triangular structure on
$H(G;G_+,G_-)$ is given by the construction above.
\end{prop}

\noindent{\em Proof}: Suppose that $R$ is a positive quasi-triangular
structure. Then $R^{-1}=(S\otimes id)R$, so that $R^{-1}$ is also
positive. Consequently, Proposition \ref{lem1} applies. In
particular, $R$ is of the form (\ref{r-form'}), and we have
properties (\ref{etaleft}) and (\ref{xiright}). Note that by
Proposition \ref{condition-construction1}, we also have properties
(\ref{etaright}) and (\ref{xileft}).

For $R$ given by (\ref{r-form'}), we have
\begin{eqnarray}
\tau\Delta\{g\}R
& = &
\sum_{h_+\in G_+}
r(h_+\,^{g_-},(g_+h_+^{-1})^{(^{h_+}g_-)}) \label{Delta-R} \\
 && 
\left\{
   h_+g_-\left(\eta((g_+h_+^{-1})^{(^{h_+}g_-)})^{(h_+\,^{g_-})}
         \right)^{-1}
\right\}\otimes
\left\{g_+h_+^{-1}(^{h_+}g_-)\xi(h_+\,^{g_-})\right\}.
 \nonumber
\end{eqnarray}
On the other hand, in the product
\[
R\Delta\{g\}
=\left(\sum_{u,v\in G_+}
    r(u,v)\{u\left(\eta(v)^{u}\right)^{-1}\}\otimes\{v\xi(u)\}\right) 
 \left(\sum_{h_+\in G_+}
    \{g_+h_+^{-1}(^{h_+}g_-)\}\otimes\{h_+g_-\}\right),
\]
we must have
\[
h_+=v^{\xi(u)},\qquad
g_+h_+^{-1}=u^{(\eta(v)^{u})^{-1}}
\stackrel{(\ref{xiright})}{=}
u^{\eta(v^{\xi(u)})^{-1}}=
u^{\eta(h_+)^{-1}}.
\]
Therefore
\begin{eqnarray*}
u & = & (g_+h_+^{-1})^{\eta(h_+)}, \\
v & = &
h_+^{\xi\left((g_+h_+^{-1})^{\eta(h_+)}\right)^{-1}}
\stackrel{(\ref{etaright})}{=}
h_+^{\left(\xi(g_+h_+^{-1})^{h_+}\right)^{-1}} \\
& \stackrel{(\ref{inverse})}{=} &
h_+^{\,^{h_+^{-1}}\left(\xi(g_+h_+^{-1})^{-1}\right)}
\stackrel{(\ref{compatible})}{=}
\left((h_+^{-1})^{\left(\xi(g_+h_+^{-1})^{-1}\right)}\right)^{-1}
\stackrel{(\ref{inverse})}{=}
\,^{\xi(g_+h_+^{-1})}h_+, \\
\eta(v)^u 
& \stackrel{(\ref{xiright})}{=} &
\eta(v^{\xi(u)})=\eta(h_+), \\
\xi(u) & = &
\xi((g_+h_+^{-1})^{\eta(h_+)})  
\stackrel{(\ref{etaright})}{=} 
\xi(g_+h_+^{-1})^{h_+},
\end{eqnarray*}
and
\begin{eqnarray}
R\Delta\{g\} & = &
\sum_{h_+\in G_+}
r((g_+h_+^{-1})^{\eta(h_+)},\,^{\xi(g_+h_+^{-1})}h_+) 
\label{R-Delta} \\ 
&& 
\left\{((g_+h_+^{-1})^{\eta(h_+)})\eta(h_+)^{-1}(^{h_+}g_-)\right\}
\otimes
\left\{(\,^{\xi(g_+h_+^{-1})}h_+)(\xi(g_+h_+^{-1})^{h_+})g_-\right\}.
 \nonumber
\end{eqnarray}

Observe that in the $G_+$-parts of each term in
$\tau\Delta\{g\}R$, we have $(g_+h_+^{-1})h_+=g_+$. Therefore by 
$\tau\Delta\{g\}R=R\Delta\{g\}$ and $r>0$, we see that for any
$g_+,h_+\in G_+$, we have
\[
\left(\,^{\xi(g_+h_+^{-1})}h_+\right)
\left((g_+h_+^{-1})^{\eta(h_+)}\right)=g_+.
\]
This is exactly (\ref{etaxi}).

We now compare $\tau\Delta\{g\}R$ and $R\Delta\{g\}$ term by term. In
order to avoid confusion, we change the index $h_+$ in $\tau\Delta\{g\}R$
to $\bar{h}_+$.

The equality $\tau\Delta\{g\}R=R\Delta\{g\}$ and $r>0$ suggests us to
consider the map
\begin{equation}\label{rel22} 
h_+\mapsto \bar{h}_+=(g_+h_+^{-1})^{\eta(h_+)}.
\end{equation}
By using (\ref{etaright}), (\ref{xileft}), and (\ref{etaxi}), we
can show that the map has the following inverse
\begin{equation}\label{rel22'} 
\bar{h}_+\mapsto h_+=(g_+\bar{h}_+^{-1})^{\xi(\bar{h}_+)}.
\end{equation}
This means that the term in $\tau\Delta\{g\}R$ indexed by $\bar{h}_+$ and
the term in $R\Delta\{g\}$ indexed by $h_+$ must be equal. By comparing
the coefficients, the $G_+$-components, and the $G_-$-components of the
corresponding terms, we have
\begin{eqnarray}
r(\bar{h}_+\,^{g_-},(g_+\bar{h}_+^{-1})^{(^{\bar{h}_+}g_-)})
& = & 
r((g_+h_+^{-1})^{\eta(h_+)},\,^{\xi(g_+h_+^{-1})}h_+)
\label{rel21} \\
g_-\left(\eta\left((g_+\bar{h}_+^{-1})^{(^{\bar{h}_+}g_-)}
            \right)^{(\bar{h}_+\,^{g_-})}
   \right)^{-1}
& = &
\eta(h_+)^{-1}(^{h_+}g_-)
\label{rel23} \\
g_+\bar{h}_+^{-1}
& = &
\,^{\xi(g_+h_+^{-1})}h_+ 
\label{rel24} \\
(^{\bar{h}_+}g_-)\xi(\bar{h}_+\,^{g_-})
& = &
(\xi(g_+h_+^{-1})^{h_+})g_-.
\label{rel25}
\end{eqnarray}
where $\bar{h}_+$ is given by (\ref{rel22}).

As pointed out earlier, (\ref{rel22}) and (\ref{rel24}) implies
(\ref{etaxi}).

Let $h_+=e$. Then from (\ref{rel22}) we have $\bar{h}_+=g_+$, so that
(\ref{rel25}) becomes 
\[
(^{g_+}g_-)\xi(g_+\,^{g_-})=\xi(g_+)g_-.
\]
This is (\ref{etaneg}).

We have from  (\ref{etaright}) and (\ref{rel22}) that
\[
\xi(\bar{h}_+)=\xi((g_+h_+^{-1})^{\eta(h_+)})
=\xi(g_+h_+^{-1})^{h_+},
\]
so that
\begin{eqnarray*}
h_+\,^{\xi(\bar{h}_+)^{-1}} 
& = &
h_+\,^{\left( \xi(g_+h_+^{-1})^{h_+} \right)^{-1}}
\stackrel{(\ref{inverse})}{=}
h_+^{\,^{h_+^{-1}}\left(\xi(g_+h_+^{-1})^{-1} \right)} 
\stackrel{(\ref{compatible})}{=} 
\left( (h_+^{-1})^{\xi(g_+h_+^{-1})^{-1}} \right)^{-1} \\
&  \stackrel{(\ref{inverse})}{=} &
\,^{\xi(g_+h_+^{-1})}h_+ 
\stackrel{(\ref{etaxi})}{=}
(g_+h_+^{-1})h_+\left((g_+h_+^{-1})^{\eta(h_+)}\right)^{-1}
\stackrel{(\ref{rel22})}{=}
g_+\bar{h}_+^{-1}.
\end{eqnarray*}
Thus $(g_+\bar{h}_+^{-1})^{\xi(\bar{h}_+)}=h_+$, and in (\ref{rel23})
we have
\begin{eqnarray*}
\eta\left((g_+\bar{h}_+^{-1})^{(^{\bar{h}_+}g_-)}
   \right)^{(\bar{h}_+\,^{g_-})}
& \stackrel{(\ref{xiright})}{=} & 
\eta\left((g_+\bar{h}_+^{-1})^{(^{\bar{h}_+}g_-)\xi(\bar{h}_+\,^{g_-})}
    \right) \\
& \stackrel{(\ref{etaneg})}{=} &
\eta\left((g_+\bar{h}_+^{-1})^{\xi(\bar{h}_+)g_-}
    \right)
=
\eta(h_+\,^{g_-}).
\end{eqnarray*}
Therefore (\ref{rel23}) becomes
\[
\eta(h_+)g_-
=
(^{h_+}g_-)\eta\left((g_+\bar{h}_+^{-1})^{(^{\bar{h}_+}g_-)}
            \right)^{(\bar{h}_+\,^{g_-})}
=
(^{h_+}g_-)\eta(h_+\,^{g_-}).
\]
This is (\ref{xineg}).

Finally, substituting (\ref{rel22}) and (\ref{rel24}) into
(\ref{rel21}) gives
\[
r(\bar{h}_+\,^{g_-},(g_+\bar{h}_+^{-1})^{(^{\bar{h}_+}g_-)})
= 
r(\bar{h}_+,g_+\bar{h}_+^{-1}).
\]
This is (\ref{rneg}).

This completes the proof that any positive quasi-triangular structure must
be given by the construction in the proposition.

Conversely, given homomorphisms $\eta$, $\xi$, and a function $r$
satisfying the conditions of the proposition, we want to show that the
formula (\ref{r-form'}) satisfies (\ref{eq_prime}), (\ref{eq_delta-R}),
and  (\ref{eq_epsilon-R}).

First of all, $R$ satisfies (\ref{eq_delta-R}) by Proposition \ref{lem1}. 

Secondly, we have
\[
(\epsilon\otimes id)R=\sum_{v\in G_+}r(e,v)\{v\}.
\]
One other hand, putting $u=e$ in (\ref{rxi}) shows that 
$r(e,?): G_+\rightarrow {\bf R}^{>0}$ is a group homomorphism. Since $G_+$
is finite, we see that $r(e,v)=1$ for all $v$. Consequently, we have
$(\epsilon\otimes id)R=1$. Similarly, we have
$(id\otimes\epsilon)R=1$. 

Finally, we may use the conditions (\ref{etaright}),
(\ref{xileft}), and (\ref{etaxi}) of the proposition to show that
(\ref{rel22}) and (\ref{rel22'}) are inverse to each other, just
as what we have done in the first part. This implies that
(\ref{rel22}) and (\ref{rel22'}) give a one-to-one  correspondence
between the terms in $\tau\Delta\{g\}R$ and
$R\Delta\{g\}$. Thus in order to show $\tau\Delta\{g\}R=R\Delta\{g\}$, it
remains to verify (\ref{rel21}-\ref{rel25}). The detailed computation is
almost the same as what we have done in the first part of the proof,
except for (\ref{rel25}). The following computation verifies (\ref{rel25}):
\[
(^{\bar{h}_+}g_-)\xi(\bar{h}_+\,^{g_-})
\stackrel{(\ref{etaneg})}{=}
\xi(\bar{h}_+)g_-
\stackrel{(\ref{rel22})}{=}
\xi((g_+h_+^{-1})^{\eta(h_+)})g_-
\stackrel{(\ref{etaright})}{=}
(\xi(g_+h_+^{-1})^{h_+})g_-.
\]
This completes the proof of the converse.

\hfill$\Box$

It remains to study the conditions (\ref{reta}), (\ref{rxi}), and
(\ref{rneg}) imposed on $r$.

\begin{prop}\label{lem4}
In a positive quasi-triangular structure
\[
R=\sum_{u,v\in G_+}r(u,v)
\{u\left(\eta(v)^{u}\right)^{-1}\}\otimes \{v\xi(u)\}
\]
on $H(G;G_+,G_-)$, we must have $r(u,v)=1$.
\end{prop}

\noindent{\em Proof}: By Proposition \ref{lem2}, the maps $\xi$
and $\eta$ in the positive quasi-triangular structure must be
homomorphisms satisfying
$(\ref{etaright}-\ref{xineg'})$. By Proposition
\ref{condition-construction1}, we conclude that we are in the
situation described in Theorem \ref{r-classification}. In
particular, $G_+'$ is a subgroup.

The one-to-one correspondence 
\[
u\mapsto u\xi(u^{-1}):\quad G_+\rightarrow G_+'
\] 
translates the group structure (\ref{rel52}) on $G_+'$ to a group
structure
\begin{equation}\label{gplus''}
u\star v = u(^{\xi(u^{-1})}v)
\end{equation}
on $G_+$. We have
\[
r(u\star w, v) 
=
r(u(^{\xi(u^{-1})}w),v)
\stackrel{(\ref{reta})}{=}
r(u,v)r(^{\xi(u^{-1})}w,v^{\xi(u)}).
\]
Let $x=\,^{w^{-1}}\xi(u)$. Then 
\[
\,^w x = \xi(u), \qquad
w^x = w^{(^{w^{-1}}\xi(u))} 
\stackrel{(\ref{compatible})}{=}
\left((w^{-1})^{\xi(u)}\right)^{-1}
\stackrel{(\ref{inverse})}{=}
\,^{\xi(u^{-1})}w.
\]
Therefore
\[
r(u\star w, v) 
=
r(u,v)r(w^x,v^{(^w x)})
\stackrel{(\ref{rneg})}{=}
r(u,v)r(w,v).
\]
Thus for fixed $v$, $r(?,v): (G_+,\star)\rightarrow {\bf R}^{>0}$
is a group homomorphism. Since $G_+$ is finite, we see that $r=1$.

\hfill$\Box$

\section{Comparing positive quasi-triangular structures}
\label{sec_gd3}

Let $G=G_+G_-$ and $G=G_+'G_-$ be two unique
factorizations of a finite group $G$. Then we have two Hopf
algebra structures $H(G;G_+,G_-)$ and
$H(G;G_+',G_-)$ on ${\bf C}G$. We show that the two
Hopf algebra structures are quasi-isomorphic.

The two unique factorizations give rise to a map
$\sigma: G_+\rightarrow G_-$ such that
\begin{equation}\label{shift}
G_+'=\{\sigma(u)u: u\in G_+\}.
\end{equation}
The fact that $G_+'$ is a subgroup implies that for
any $u,v\in G_+$,
\begin{equation}\label{group-condition}
\sigma((u^{\sigma(v)})v)=\sigma(u)(^u\sigma(v)).
\end{equation}

To avoid confusion about two structures on the same space 
${\bf C}G$. We consider $G=G_+G_-$ as the ``standard''
factorization, and $G=G_+'G_-$ as the ``shifting'' of
the standard factorization. All the notations in Section
\ref{sec_gd} refer to operations relative to the unique
factorization $G=G_+G_-$ and the Hopf algebra structure
$H(G;G_+,G_-)$. For any $g\in G$, we use $\{g\}$ and
$\{g\}'$ to denote $g$ considered as an element in
$H(G;G_+,G_-)$ and in $H(G;G_+',G_-)$, respectively.

By solving
\[
\sigma(u)ux=y\sigma(v)v,\qquad
u,v\in G_+,\quad x,y\in G_-,
\]
for $v$ and $y$, we have 
\[
v=u^x,\quad y=\sigma(u)(^ux)\sigma(v)^{-1}=\sigma(u)(^ux)\sigma(u^x)^{-1}.
\]
Therefore the left action of $G_+'$ on $G_-$ and the
right action of $G_-$ on $G_+'$ are given by 
\begin{eqnarray}
G_+'\times G_-\rightarrow G_+',  && 
(\sigma(u)u,x)\mapsto \sigma(u^x)u^x, \label{action'1}\\
G_+'\times G_-\rightarrow G_-,  && 
(\sigma(u)u,x)\mapsto \sigma(u)(^ux)\sigma(u^x)^{-1}. 
\label{action'2}
\end{eqnarray}

\begin{prop}\label{quasi-isomorphism}
Denote
\begin{equation}\label{quasi-iso}
\phi\{\sigma(u)ux\}'=\{ux\}:\quad 
H(G;G_+',G_-)\rightarrow H(G;G_+,G_-),
\end{equation}
and
\begin{equation}\label{twist}
T=\sum_{u,v\in G_+}\{u\sigma(v)\}\otimes\{v\}.
\end{equation}
Then $(\phi,T)$ is a quasi-isomorphism of Hopf algebras.
\end{prop}

\noindent{\em Proof}: We first show that $\phi$ is an isomorphism of
algebras. From the action (\ref{action'1}), we see that
\[
\{\sigma(u)ux\}'\{\sigma(v)vy\}'=
\delta_{u^x,v}\{\sigma(u)uxy\}'.
\]
It follows from this that $\phi$ preserves the multiplication. It is also
easy to see that $\phi$ preserves the unit.

Next we show that 
\begin{equation}\label{quasi-isomorphism2}
(\phi\otimes\phi)\Delta(\{g\}')=
T(\Delta\phi(\{g\}'))T^{-1}. 
\end{equation}
It is easy to verify that 
\[
T^{-1}=\sum_{u,v\in G_+}\{u\sigma(v)^{-1}\}\otimes\{v\}.
\]
Then for any $g_+\in G_+,g_-\in G_-$, we have
\[
T(\Delta\phi\{\sigma(g_+)g_+g_-\}')T^{-1} =
\sum_{h_+\in G_+}\left\{\left((g_+h_+)^{\sigma(h_+)^{-1}}\right)
     \sigma(h_+)(^{h_+}g_-)\sigma(h_+^{g_-})^{-1}\right\}\otimes \{h_+g_-\}.
\]
On the other hand, 
\[
\Delta\{\sigma(g_+)g_+g_-\}'   
\stackrel{(\ref{action'2})}{=} 
\sum_{h_+\in G_+}\{\sigma(g_+)g_+(\sigma(h_+)h_+)^{-1}
  \sigma(h_+)(^{h_+}g_-)\sigma(h_+^{g_-})^{-1}\}'
  \otimes\{\sigma(h_+)h_+g_-\}'.
\]
Since $G_+'$ is a subgroup, we have
\[
\sigma(g_+)g_+(\sigma(h_+)h_+)^{-1}=\sigma(w)w
\]
for some $w\in G_+$. By considering the $G_+$-components in the unique
factorization $G=G_-G_+$, we find $w=(g_+h_+^{-1})^{\sigma(h_+)^{-1}}$.
Therefore 
\[
\sigma(g_+)g_+(\sigma(h_+)h_+)^{-1}
=\sigma\left((g_+h_+^{-1})^{\sigma(h_+)^{-1}}\right)
 (g_+h_+^{-1})^{\sigma(h_+)^{-1}},
\]
and
\[
(\phi\otimes\phi)\Delta\{\sigma(g_+)g_+g_-\}' 
=
\sum_{h_+\in G_+}
\left\{\left((g_+h_+^{-1})^{\sigma(h_+)^{-1}}\right)
  \sigma(h_+)(^{h_+}g_-)\sigma(h_+^{g_-})^{-1}\right\}
\otimes\{h_+g_-\}.
\]
This completes the verification of (\ref{quasi-isomorphism2}).

Finally, it is easy to compute the following
\begin{eqnarray*}
(T\otimes 1)(\Delta\otimes id)T & = &
\sum_{u,v,w\in G_+}
\{u\sigma(v)(^v\sigma(w))\}\otimes
\{v\sigma(w)\}\otimes
\{w\}, \\
(1\otimes T)(id\otimes\Delta)T & = &
\sum_{u,v,w\in G_+}
\{u\sigma(v^{\sigma(w)}w)\}\otimes
\{v\sigma(w)\}\otimes
\{w\}.
\end{eqnarray*}
It then follows from (\ref{group-condition}) that 
$(T\otimes 1)(\Delta\otimes id)T
=(1\otimes T)(id\otimes\Delta)T$, so that the compatibility
condition is verified.

\hfill$\Box$

Now we apply the proposition above to the special case in
Theorem \ref{r-classification}, with $\sigma(u)=\xi(u^{-1})$.
Note that by taking inverse, $G_+'$ is the same as the one
given in Theorem \ref{r-classification}. Thus the
quasi-isomorphism constructed in Proposition
\ref{quasi-isomorphism} translates the quasi-triangular
structure $R$ on $H(G;G_+,G_-)$ into another quasi-triangular
structure 
\[
R'=(\phi\otimes\phi)^{-1}((\tau T)RT^{-1})
\]
on $H(G;G_+',G_-)$. An easy computation gives
\[
(\tau T)RT^{-1}=
\sum_{u,v\in G_+}\{u\eta(v)^{-1}\xi(v)\}\otimes\{v\},
\]
so that
\begin{equation}\label{r'}
R'=
\sum_{u,v\in G_+}\{\xi(u^{-1})u\eta(v)^{-1}\xi(v)\}'
\otimes\{\xi(v^{-1})v\}'.
\end{equation}
By Theorem \ref{r-classification}, $R'$ is given by
homomorphisms $\xi',\eta':G_+'\rightarrow G_-$. From
the second component in (\ref{r'}), we have
$\xi'(\xi(v^{-1})v)=e$. By (\ref{xiright}), the triviality of
$\xi'$ implies that the first component in (\ref{r'}) is of the
form
$\{\xi(u^{-1})u\eta'(\xi(v^{-1})v)^{-1}\}'$. Therefore we
conclude that 
\begin{equation}\label{normal}
\xi'(\xi(u^{-1})u)=e,\qquad
\eta'(\xi(v^{-1})v)=\xi(v)^{-1}\eta(v).
\end{equation}

\bigskip

\noindent{\bf Definition} Let $R$  be a positive
quasi-triangular structure on $H(G; G_+ , G_-)$ given by
homomorphisms
$\xi,\eta: G_+ \to G_- $ as in Theorem \ref{r-classification}. 
We say that $R$ is {\em normal} if $\xi(u) = e $ for all 
$u\in G_+$. We say that the pair $(H(G;G_+,G_-), R)$ is {\em
normal} if $R$ is normal. 
       
\bigskip

Thus in the special case described in Theorem \ref{r-classification},
Proposition \ref{quasi-isomorphism} implies the following.

\begin{thm}\label{quasi-normalize}
Every pair $(H, R)$, where $H$ is a finite
dimensional Hopf algebra with a positive basis and
$R$ is a positive quasi-triangular structure in this basis, 
is 
quasi-isomorphic to a normal one.
\end{thm}

\hfill$\Box$

\section{Positive triangular structures}

A quasi-triangular structure $R$ is {\em triangular} if it
further satisfies $(\tau R)R=1\otimes 1$. For the positive
quasi-triangular structure $R$ given by Theorem
\ref{r-classification}, we have 
\[
(\tau R)R=\sum
\{v\xi(u)\left(\eta(\bar{v})^{\bar{u}}\right)^{-1}\}\otimes
\{u\left(\eta(v)^{u}\right)^{-1}\xi(\bar{u})\},
\]
where the summation is over all $u,v,\bar{u},\bar{v}\in G_+$
satisfying
\begin{equation}\label{tri1}
v^{\xi(u)}=\bar{u},\qquad \,^{\eta(v)}u=\bar{v}.
\end{equation}
Since $1\otimes 1=\sum_{u,v\in G_+}\{u\}\otimes\{v\}$, we see
that 
$R$ is triangular if and only if (\ref{tri1}) implies
\begin{equation}\label{tri2}
\xi(u)=\eta(\bar{v})^{\bar{u}},\qquad
\eta(v)^{u}=\xi(\bar{u}).
\end{equation}
Note that the first equality in (\ref{tri1}) implies 
$\eta(v)^{u}=\eta(v^{\xi(u)})=\eta(\bar{u})$. Therefore under
the assumption (\ref{tri1}), the second equality of
(\ref{tri2}) is equivalent to $\xi=\eta$. Furthermore, the
property $\xi=\eta$ and (\ref{tri1}) imply
\[
\eta(\bar{v})\bar{u}=\eta(^{\eta(v)}u)\bar{u}
=\xi(^{\eta(v)}u)\bar{u}\stackrel{(\ref{etaleft})}{=}
(^v\xi(u))(v^{\xi(u)})=v\xi(u).
\]
In particular, the first equality of (\ref{tri2})
also holds. Thus we conclude that $R$ is triangular if and
only if $\xi=\eta$.

\begin{thm}\label{tri}
There is a one-to-one correspondence between the following data
\begin{enumerate}
\item a finite dimensional Hopf algebras with a  positive basis
and a positive triangular structure;
\item a unique factorization $G=G_+G_-$ of finite group, and a
homomorphism $\xi: G_+\rightarrow G_-$ such that
$A=\{u\xi(u^{-1}):u\in G_+\}$ is an abelian normal subgroup;
\item a unique factorization $G=G_+G_-$ of finite group, and a
homomorphism $\xi: G_+\rightarrow G_-$ satisfying 
$uv=(^{\xi(u)}v)(u^{\xi(v)})$ and 
$\xi(^x u)x\,^u=x\xi(u)$;
\item a group $G_-$, an abelian group $A$ acted upon by $G_-$
as automorphisms, and a 1-cycle $\zeta$ of $G_-$ with
coefficient in $A$.
\end{enumerate}
\end{thm}

\noindent{\em Proof}: We already know the first item is
equivalent to a unique factorization $G=G_+G_-$ of a finite
group with a homomorphism $\xi: G_+\rightarrow G_-$ such
that 
\[
G_+^{'} \, = \, \{u\xi(u^{-1}): \, \, u \in  G_+\},
\qquad
G_+^{''} \, = \, \{\xi(u^{-1})u: \, \, u\in G_+\}
\]
are normal subgroups, and
\[
F(u\xi(u^{-1}))=\xi(u)u^{-1}: \quad G_+'\rightarrow G_+''
\]
is a homomorphism. Since $\xi(u)u^{-1}=(u\xi(u^{-1}))^{-1}$,
we see that $G_+^{'}=G_+^{''}$ and $F(a)=a^{-1}$. Therefore
the condition is equivalent to $A=G_+^{'}=G_+^{''}$ being an
abelian normal subgroup. This proves the equivalence between
the first two items.

By Proposition \ref{condition-construction1}, if $\xi$
induces a triangular structure, then the two conditions in the
third item must be satisfied. Conversely, we need to show
that the two conditions imply the other conditions in 
Proposition \ref{condition-construction1}. For example, the
two conditions imply
\[
\xi(^{\xi(u)}v)\xi(u)^v=\xi(u)\xi(v)=\xi(uv)=
\xi(^{\xi(u)}v)\xi(u^{\xi(v)}).
\]
Canceling the left factor, we get (\ref{etaright}). The other
conditions can be verified similarly. This proves the
equivalence to the third item.

Finally, the fourth item is the reinterpretation in terms of
the alternative description. From the discussion above, we see
that a positive triangular structure means $A=G_+'$ is
abelian, and
\[
F(ax)=a^{-1}x,\qquad a\in A, \, x\in G_-.
\]
Since these already implies the conditions in the second item
are satisfied, we see that there is no further condition on the
1-cycle $\zeta$. This proves the equivalence to the fourth
item.

\hfill$\Box$

Now let us apply the theory of Section \ref{sec_gd3} to
normalize positive triangular structures. Since $\xi=\eta$, 
both $\xi'$ and $\eta'$ are trivial by
(\ref{normal}). Another way to see the triviality is to use
the fact that a quasi-isomorphism carries triangular
structures to triangular structures. Therefore
$(\xi',\eta')$ must induce a triangular
structure, which implies $\xi'=\eta'$. Since $\xi'$ is
trivial, so is $\eta'$. Thus, we find $R'=1'\otimes 1'$, and
the triangular Hopf algebra $(H(G;G_+,G_-),R)$ is isomorphic to
the twisting of the triangular Hopf algebra
$(H(G;A,G_-),1'\otimes 1')$, with the twist given by
\[
T'=(\phi\otimes\phi)^{-1}(T)
=\sum_{u,v\in G_+}\{\xi(u^{-1})u\xi(v^{-1})\}'
   \otimes\{\xi(v^{-1})v\}'
=\sum_{a,b\in A}\{a\zeta(b^{-1})\}'
   \otimes\{b\}'.
\]
Finally, we note that since $R'=1'\otimes 1'$, $H(G;A,G_-)$ is
cocommutative. Thus $H(G;A,G_-)$ is a group algebra, and we
conclude that any positive triangular structure is the
twisting of a group algebra. Explicit formulae for this
conclusion can be found in Section 4 of \cite{eg}.

\section{Groupoids and the set-theoretical Yang-Baxter
equation}
\label{set-theory1}

In \cite{lyz:positive}, we have shown that the positivity
condition on a Hopf algebra implies that the Hopf algebra is
essentially set-theoretical. Theorem \ref{r-classification}
says that positive quasi-triangular structures on such Hopf
algebras are also set-theoretical. As a result, we expect our
theory to lead to set-theoretical solutions of the Yang-Baxter
equation.  In this section, we explain how this happens.

Recall \cite{mk:book} that a {\em groupoid} over a set $B$ (called {\em base
space}) is a set $\Gamma$ (called {\em total space}) together with 
\begin{enumerate}
\item two surjections $\alpha,\beta: \Gamma\rightarrow B$;
\item a product $\mu: (\gamma_1,\gamma_2)\mapsto \gamma_1\gamma_2$ in
$\Gamma$, defined when $\beta(\gamma_1)=\alpha(\gamma_2)$;
\item an identity map $e: b\mapsto e_b, B\rightarrow\Gamma$;
\item an inversion map $\sigma: \gamma\mapsto\gamma^{-1}, 
\Gamma\rightarrow\Gamma$
\end{enumerate}
such that
\[
\alpha(\gamma_1\gamma_2)=\alpha(\gamma_1),\quad
\beta(\gamma_1\gamma_2)=\beta(\gamma_2),\quad
\alpha(e_b)=\beta(e_b)=b,
\]
\[
\alpha(\gamma^{-1})=\beta(\gamma),\quad
\beta(\gamma^{-1})=\alpha(\gamma),
\]
and the usual axioms similar to those for groups are satisfied. If $\Gamma$
is finite, then we have an algebra structure on the vector space
${\bf C}\Gamma$ 
\[
\{\gamma_1\}\{\gamma_2\} = 
\left\{\begin{array}{ll}
\{\gamma_1 \gamma_2\} & \mbox{if $\beta(\gamma_1) = \alpha(\gamma_2)$} \\
0 & \mbox{if $\beta(\gamma_1) \neq \alpha(\gamma_2)$}
\end{array}
\right., \qquad
e= \sum_{b \in B} e_b,
\]
called the {\em linearization} of the groupoid $\Gamma$ or the 
{\em groupoid algebra} of $\Gamma$.

For example, given a unique factorization $G=G_+G_-$, we have
the following groupoid
$\Gamma_+$ with $G$ as the total space and with $G_+$ as the base space:
\begin{enumerate}
\item $\alpha_{+}$: $g\mapsto g_+,\, G \rightarrow G_+$, and
      $\beta_{+}$: $g\mapsto\bar{g}_+,\, G \rightarrow G_+$;
\item $\mu_+: (g,h)\mapsto gh_-$ when $\bar{g}_+=h_+$;
\item $e_{+}: g_+\mapsto g_+,\, G_+\rightarrow G$;
\item $\sigma_+: g\mapsto\bar{g}_+g_-^{-1}=\bar{g}_-^{-1}g_+,
       \, G \rightarrow G$.
\end{enumerate}
The linearization of $\Gamma_+$ is the algebra structure of $H(G;G_+,G_-)$.

We can read set-theoretical
information about the groupoid from its linearization.
Specifically, let $\Gamma$ be a (finite) groupoid over $B$.
For a  {\em positive element}
\[
a=\sum_{\gamma\in\Gamma}r(\gamma)\{\gamma\},\qquad
r(\gamma)\geq 0
\]
of the groupoid algebra, we denote
\[
L(a)=\{\gamma:\quad r(\gamma)>0\}.
\]
Since $\Gamma$ is a positive basis of the groupoid algebra,
$L$ has the following property: If $a_1$ and $a_2$ are two
positive elements, then
\begin{equation}\label{mul0}
L(a_1a_2)=\{\gamma_1\gamma_2:\quad 
  \gamma_1\in L(a_1), \,
  \gamma_2\in L(a_2), \,
  \beta(\gamma_1)=\alpha(\gamma_2)\}.
\end{equation}
If we define the product of two subsets
$L_1,L_2\sub\Gamma$ as 
\begin{equation}\label{mul}
L_1L_2=\{\gamma_1\gamma_2:\quad 
  \gamma_1\in L_1, \,
  \gamma_2\in L_2, \,
  \beta(\gamma_1)=\alpha(\gamma_2)\},
\end{equation}
then (\ref{mul0}) becomes $L(a_1a_2)=L(a_1)L(a_2)$.

The subset
\[
E_{\Gamma}=\{e_b:\quad b\in B\}\subset \Gamma
\]
is clearly a unit of the multiplication (\ref{mul}). The next result tells
us which subset of $\Gamma$ is invertible.

\begin{prop}\label{bisect}
Let $L\sub\Gamma$ be a subset of a groupoid  $\Gamma$ over $B$. Then the
following are equivalent
\begin{enumerate}
\item There is a subset $K\sub\Gamma$ such that $LK=E_{\Gamma}$ and
$KL=E_{\Gamma}$;
\item The restrictions $\alpha|_L,\beta|_L: L\rightarrow B$ are
bijections.
\end{enumerate}
\end{prop}

\noindent{\em Proof}: 
We show that the first statement implies the second.

Since $LK=E_{\Gamma}$, for any $b\in B$,
there are $\gamma_1\in L$ and $\gamma_2\in K$ such that
$\beta(\gamma_1)=\alpha(\gamma_2)$ and $\gamma_1\gamma_2=e_b$. In
particular, we have $\alpha(\gamma_1)=\alpha(\gamma_1\gamma_2)=b$ and
$\beta(\gamma_2)=\beta(\gamma_1\gamma_2)=b$. Thus $\alpha|_L$
and $\beta|_K$ are surjective.

Now suppose we have $\gamma_1,\gamma_1'\in L$ such that
$\alpha(\gamma_1)=\alpha(\gamma_1')=a$. By the surjectivity of $\beta|_K$, 
we can find $\gamma_2\in K$ such that $\beta(\gamma_2)=a$. Since
$KL=E_{\Gamma}$, we conclude that $\gamma_2\gamma_1=e_b$ and
$\gamma_2\gamma_1'=e_{b'}$ for some $b,b'\in B$. Then
\[
\beta(\gamma_1)=\beta(\gamma_2\gamma_1)=b
=\alpha(\gamma_2\gamma_1)=\alpha(\gamma_2)
=\alpha(\gamma_2\gamma_1')=b'
=\beta(\gamma_2\gamma_1')=\beta(\gamma_1').
\]
In particular, the products $\gamma_1\gamma_2$ and
$\gamma_1'\gamma_2$ make sense. Since
$\gamma_1\gamma_2,\gamma_1'\gamma_2\in LK=E_{\Gamma}$, we see
that $\gamma_1\gamma_2=\gamma_1'\gamma_2=e_a$. Combining this
with $\gamma_2\gamma_1=e_b=e_{b'}=\gamma_2\gamma_1'$, we
conclude that $\gamma_1=\gamma_2^{-1}=\gamma_1'$. This proves
the injectivity of $\alpha|_L$.

The bijectivity of $\beta|_L: L \rightarrow B$ can be proved
similarly.

Conversely, given the second statement, it is easy
to verify that $K=\{\gamma^{-1}: \gamma\in L\}$
satisfies the first statement.

\hfill$\Box$

A subset $L\sub\Gamma$ of a groupoid satisfying the equivalent
conditions of the proposition above is called a 
{\em bisection}. With product (\ref{mul}), the collection 
${\cal U}(\Gamma)$ of all bisections of a groupoid $\Gamma$
form a group.

As an application of Proposition \ref{bisect}, let us prove Proposition
\ref{positiveR}. We note that the tensor algebra $H(G;G_+,G_-)\otimes
H(G;G_+,G_-)$ is the linearization of the product groupoid
$\Gamma_+\times\Gamma_+$. Under the assumption of Proposition
\ref{positiveR}, both $R$ and $R^{-1}$ are positive elements. Then
$RR^{-1}=1\otimes 1=R^{-1}R$ implies that 
$L(R)L(R^{-1})=L(1\otimes 1)=
E_{\Gamma_+\times\Gamma_+}=L(R^{-1})L(R)$, i.e., $L(R)$ is
invertible. By Proposition \ref{bisect}, the restriction of
$\alpha_{\Gamma_+\times\Gamma_+}(g,g)=(g_+,g_+)$ on ${\cal R}=L(R)$ is
bijective.

As another application, we consider the positive
quasi-triangular structure $R$ in the classification Theorem
\ref{r-classification}. Denote
\begin{equation}\label{r-mother}
{\cal R}=L(R)=
\left\{\left(u\left(\eta(v)^{u}\right)^{-1},\,v\xi(u)\right):
\quad u, v \in G_+\right\}.
\end{equation}
Then we clearly have
\[
L(R_{12})={\cal R}\times\{e\}={\cal R}_{12},\qquad\mbox{etc.}
\]
Applying $L$ to the Yang-Baxter equation satisfied by the positive
quasi-triangular structure, we see that ${\cal R}$ also satisfies the
following {groupoid-theoretical Yang-Baxter equation} introduced in
\cite{w-x:R}
\begin{equation}\label{ybe-set}
{\cal R}_{12} {\cal R}_{13} {\cal R}_{23} \, = \, {\cal R}_{23} {\cal R}_{13} {\cal R}_{12} ,
\end{equation}
which is an equality inside the group 
${\cal U}(\Gamma_+\times\Gamma_+\times\Gamma_+)$ of bisections.

To get set-theoretical solutions of the Yang-Baxter equation
from (\ref{ybe-set}), we recall that the quasi-triangular
structure $R$ induces a solution of the Yang-Baxter equation
on any $H$-module. The set-theoretical analogue of modules
is sets acted upon by groupoids.

Let $\Gamma$ be a groupoid over $B$. A (left) {\em $\Gamma$-set} consists of
a set $X$, a map $J:X\rightarrow B$, and an action
\[
(\gamma,x)\mapsto \gamma x \in X,
\qquad
\mbox{for $\gamma\in\Gamma$, $x\in X$, satisfying 
$\beta(\gamma)=J(x)$}.
\]
The action is required to satisfy
\[
(\gamma_1\gamma_2) x = \gamma_1(\gamma_2 x),
\qquad
e_{J(x)}x=x,
\qquad
J(\gamma x)=\alpha(\gamma),
\]
whenever the relevant actions are defined. The vector space ${\bf C}X$ has
an obvious module structure over the linearization of $\Gamma$.

For any $L\in{\cal U}(\Gamma)$ and $x\in X$, the equation
\[
Lx=\gamma x,\qquad
\mbox{for the unique $\gamma\in L$ satisfying $\beta(\gamma)=J(x)$}
\]
defines a (left) action of the group ${\cal U}(\Gamma)$ of bisections on the
set $X$. Now if ${\cal R}\in{\cal U}(\Gamma\times\Gamma)$ satisfies the
groupoid-theoretical Yang-Baxter equation (\ref{ybe-set}), then the map
${\cal R}_X: X\times X\rightarrow X\times X$ induced by ${\cal R}$ is a
set-theoretical solution of the Yang-Baxter equation over $X$.

Now we compute the set-theoretical solution of the Yang-Baxter
equation induced by the action of the bisection
(\ref{r-mother}) on the simplest $\Gamma_+$-set, the unit
$\Gamma_+$-set $id: G_+\rightarrow G_+$. In this case,
${\cal R}_{G_+}$ is given by the following diagram
\[
\begin{array}{rcl}
& (u\left(\eta(v)^{u}\right)^{-1},\,v\xi(u)) & \\
{\scriptstyle \beta_+\times\beta_+}\swarrow &&
\searrow {\scriptstyle \alpha_+\times\alpha_+} \\
(u^{\left(\eta(v)^{u}\right)^{-1}},v^{\xi(u)}) &
\stackrel{\scriptstyle {\cal R}_{G_+}}{---\longrightarrow} &
(u,v)
\end{array}
\]
By $u^{\left(\eta(v)^{u}\right)^{-1}}=\,^{\eta(v)}u$, we have
\[
{\cal R}_{G_+}^{-1}(u,v)=(^{\eta(v)}u,v^{\xi(u)}).
\]
Solving the equation, we get the set-theoretical solution 
\begin{equation}\label{r-mtx}
{\cal R}_{G_+}(u,v)=(u^{\eta(v)},\,^{\xi(u)}v)
\end{equation}
of the Yang-Baxter equation over $G_+$.

Direct computation shows that if $\xi$ and $\eta$ are two group
homomorphisms satisfying (\ref{etaxi}), then (\ref{r-mtx}) is
already a set-theoretical solution of the Yang-Baxter
equation. Moreover, in order for (\ref{etaxi}) and
(\ref{r-mtx}) to make sense, we do not even need to know
anything about $G_-$. The only data we need are actions 
$(u,v)\mapsto \,^{\xi(u)}v$ and $(u,v)\mapsto u^{\eta(v)}$ of
$G_+$ on itself. This is the basic data for the
construction in \cite{lyz:setybe}.

In \cite{lyz:setybe}, we have an alternative description of
our set-theoretical solution of the Yang-Baxter equation in
terms of bijective 1-cocycles. If the solution is given by
(\ref{r-mtx}), then we may take the groups $G$ and $A$ in
\cite{lyz:setybe} to be $G_+$ and $G_+'$. Moreover, we may take
the action of $G$ and $A$ in \cite{lyz:setybe} to be the
translation of the action $(u,v)\mapsto \,^{\xi(u)}v$ of $G_+$
on itself under the bijection
\[
\pi(u)=u\xi(u^{-1}):\quad G_+\rightarrow G_+'.
\]
Finally, $\pi$ also serves as the bijective 1-cocycle for the 
alternative description in \cite{lyz:setybe}. This 1-cocycle
is related to the 1-cycle $\zeta$ by
\[
\zeta=\xi\circ\pi^{-1}.
\]

We would like to end the section by mentioning that a comprehensive theory
can be established for the Yang-Baxter equation on groupoids. Indeed, we
can formulate the definitions of Hopf groupoids and quasi-triangular
structures on them. 
We can further show that any unique factorization of a group
induces a Hopf groupoid, and that all of its
quasi-triangular structures are given by the bisections of the
form (\ref{r-mother}). Moreover, we can introduce the notion
of quasi-isomorphisms of Hopf groupoids and establish the
set-theoretical analogue of Theorem \ref{quasi-normalize}. In
particular, we conclude that actions of the braid group
$B_n$ on the set $X^n$ induced by the solution ${\cal R}_X$ is
equivalent to the action induced by a {\em normal} solution.
Note that by (\ref{etaxi}), the special solution (\ref{r-mtx})
given by a normal quasi-triangular structure is of the form
\[
{\cal R}_{G_+'}(u,v)=(v^{-1}uv,v),
\]
which is the conjugate solution (with respect to the new
multiplication (\ref{gplus''}) coming from $G_+'$). Thus we
conclude that the braid group action induced by any
set-theoretical solution of the Yang-Baxter equation of the
form (\ref{r-mtx}) is equivalent to the action induced from a
conjugate solution. All these considerations motivate our
paper \cite{lyz:setybe} on set-theoretical solutions of the 
Yang-Baxter equation.

\section{An example}

Any unique factorization $G=G_+G_-$
induces another unique factorization 
$\tilde{G}=G \times G=\tilde{G}_+\tilde{G}_-$, with
\[
\tilde{G}_{+}=\{(g_+,g_-): \, g_+\in G_+, g_-\in G_-\},
\qquad
\tilde{G}_{-}=\{(g, g): \, g \in G\}.
\]
The Hopf algebra induced by this unique factorization is in fact the
Drinfel'd double of $H(G;G_+,G_-)$ (see \cite{lyz:positive}).

Consider homomorphisms
\[
\left\{\begin{array}{l}
\xi(g_+,g_-)=(g_-,g_-) \\
\eta(g_+,g_-)=(g_+,g_+)
\end{array}\right. : \qquad
\tilde{G}_{+} \rightarrow \tilde{G}_{-}.
\]
The induced subgroups (as in Theorem \ref{r-classification})
$\tilde{G}_+'=G\times \{e\}$ and 
$\tilde{G}_+''=\{e\}\times G$ are clearly normal. It
is also easy to see that the map
$F:\tilde{G}_{+}'\rightarrow\tilde{G}_{+}''$ is given by
$F(a,e)=(e,a)$, which is clearly a group isomorphism. From
this we get the standard quasi-triangular structure on the Drinfel'd
double of $H(G;G_+,G_-)$.

To find the alternative description, we use the identification
\begin{equation}\label{id}
\tilde{G}_+'\cong G: \,(a,e)\leftrightarrow a;\qquad
\tilde{G}_-\cong G: \,(g,g)\leftrightarrow g.
\end{equation}
Then $(g,g)(a,e)(g,g)^{-1}=(gag^{-1},e)$ implies that $G$
acts on $A=G$ by conjugations. Since $(a,e)=u\xi(u^{-1})$ for
$u=(a_+,a_-^{-1})$, the 1-cycle is (after the identification
(\ref{id}))
\[
\zeta(a)=a_-^{-1}.
\]
Moreover, since
$(e,a)=(a^{-1},e)(a,a)$ with respect to the unique
factorization
$\tilde{G}=\tilde{G}_+'\tilde{G}_-$, the automorphism on 
$G\lcpr_{\mbox{\scriptsize conj}} G$ is
\[
F(a\lcpr g)=F(a\lcpr e)F(e\lcpr g)=
(a^{-1}\lcpr e)(e\lcpr a)(e\lcpr g)=a^{-1}\lcpr ag.
\]

The quasi-triangular structure induces a solution (\ref{r-mtx})
of the Yang-Baxter equation on the set $\tilde{G}_+$.
A detailed computation shows that this solution is the
one given in Theorem 9.2 of \cite{w-x:R} (after the
identification 
$\tilde{G}_+\cong G: (g_+,g_-)\leftrightarrow g_+g_-^{-1}$).

To find the bijective 1-cocycle description of this solution,
we use the construction near the end of last section. Thus we
take the group $G$ and $A$ in \cite{lyz:setybe} to be
\[
\tilde{G}_+=G_+\times G_-,\qquad
G\cong G\times\{e\}.
\]
The bijective 1-cocycle is
\[
\pi(g_+,g_-)=(g_+,g_-)\xi(g_+,g_-)^{-1}=
(g_+g_-^{-1},e)\in G\times\{e\}
\leftrightarrow g_+g_-^{-1} \in G.
\]
Moreover, the action of $\tilde{G}_+$ on $G$ may be computed as
follows: If $\xi(g_+,g_-)(a_+,a_-)=(b_+,b_-)(h,h)$,
then $(g_+,g_-)\cdot \pi(a_+,a_-)=\pi(b_+,b_-)$. Since
$g_-a_+=b_+h$ and $g_-a_-=b_-h$, we have
\[
(g_+,g_-)\cdot (a_+a_-^{-1})=
(b_+h)(b_-h)^{-1}=
(g_-a_+)(g_-a_-)^{-1}=
g_-(a_+a_-^{-1})g_-^{-1}.
\]
Therefore the action is given by a conjugation 
\[
(g_+,g_-)\cdot a=g_-ag_-^{-1}.
\]
This interpretation of Weinstein and Xu's solution of
the Yang-Baxter appeared in \cite{lyz:setybe}.

\end{document}